\documentclass[12pt]{article}
\usepackage{amssymb}
\begin{document}
\title{{Stochastic differential equations on noncompact manifolds:
 moment stability and its topological consequences}}
\author{Xue-Mei Li\thanks{Research supported by SERC grant GR/H67263.}\\
Mathematics Institute, University of Warwick, Coventry CV4 7AL\\
e-mail: xl@maths.warwick.ac.uk.bitnet}
\date{}
\maketitle

\newcommand{\A}{{\bf \cal A}}
\newcommand{\B}{{ \bf \cal B }}
\newcommand{\C}{{\cal C}}
\newcommand{\E}{{\Bbb E}}
\newcommand{\F}{{\cal F}}
\newcommand{\G}{{\cal G}}
\newcommand{\h}{{\cal H}}
\newcommand{\K}{{\cal K}}
\newcommand{\R}{{\Bbb R}}
\newcommand{\half}{{  {1\over 2}  }}
\newcommand{\heatsemif}{{ {\rm e}^{ \half t\Delta^{h,1}}   }}
\newcommand{\heatsemi}{{ {\rm e}^{\half t \Delta^{h}}   }}

\newtheorem{theorem}{Theorem}[section]
\newtheorem{proposition}[theorem]{Proposition}
\newtheorem{lemma}[theorem]{Lemma}
\newtheorem{corollary}[theorem]{Corollary}
\newtheorem{definition}{Definition}[section]

\def\limsup{\mathop{\overline{\rm lim}}}
\def\liminf{\mathop{\underline{\rm lim}}}
\def\exp{{\rm e}}

\noindent
{\bf Summary.}
In this paper we discuss the stability of stochastic differential equations
 and the interplay between  the moment stability of a SDE and the topology 
of the underlying manifold. Sufficient and necessary conditions
are given  for the  moment stability of a SDE in terms of the coefficients. 
Finally  we prove a  vanishing result for the  fundamental group of a complete
Riemannian manifold in terms of purely geometrical quantities.


\noindent
{\bf Running Title:} Moment Stability and its  Topological Consequences.

\section{Introduction}

The aim of this paper is to give conditions under which a stochastic dynamical
system is moment stable and to relate this to the topological properties
of the underlying space, following an approach of Elworthy \cite{ELflow}.
 To be more precise we need the following set up.

\bigskip

\noindent {\bf A.}
Let $M$ be a smooth manifold. Consider on $M$ the stochastic differential
equation (SDE) on $M$:

\begin{equation}
dx_t=X(x_t)\circ dB_t +A(x_t)dt.
\end{equation}
 Here  $B_t$ is a $m$-dimensional Brownian motion 
 on a filtered probability  space  $\{\Omega,{\cal F}, \F_t, P\}$,
 $X$ is $C^3$   from $ \R^m\times M$ to the tangent bundle $TM$
  with $X(x)$:  $\R^m$ $\to T_xM$ a  linear map for each $x$ in $M$,
  and $A$ is  a $C^2$ vector field on $M$. 

For each $x$ in $M$ there is a solution $\{F_t(x)\}$ to (1) starting from
 $x$  (maybe with  explosion time $\xi(x)$). If $\xi(x)=\infty$ a.s. for each
 $x$, we say the SDE is complete or has no explosion.  There are also the
 standard  transition semigroup $P_t$ and its infinitesimal generator $\A$.

Formally the derivative $v_t=TF_t(v_0)$ of $F_t$ at $x_0$ in the 
direction $v_0$ satisfies the derivative stochastic differential
equation on $TM$:

\begin{equation}\label{2}
dv_t=\delta X(v_t)\circ dB_t +\delta A(v_t)dt,
\end{equation}
where  $\delta X$ and $\delta A$ are obtained by differentiate $X$ and $A$
respectively (with a twist). See \cite{ELbook}. Furthermore $v_t$  has the
 same explosion time as $x_t$.

If $M$ is given a Riemannian structure with Levi-Civita connection, then
(\ref{2}) is  equivalent to  (1) together with the covariant equation
along the paths of $\{x_t\equiv F_t(x_0): t\ge 0\}$ 
\begin{equation}\label{3}
Dv_t=\nabla X(v_t)\circ dB_t +\nabla A(v_t)dt.
\end{equation}

 If the SDE is strongly complete, i.e.  if there is a version of 
 $\{F_t(x)\}$ which is jointly continuous in $t$ and
 $x$ for almost all $\omega \in \Omega$, then $TF_t$  is the  derivative of
 $F_t(x)$ in the classical sense\cite{ELbook}.  But this is not  needed here. 
 What we  really need  is the  concept of  strong 1-completeness:

A SDE is {\it strongly 1-complete} if for each compact smooth curve in $M$,
 $\{F_t(x)\}$ has a version which is jointly continuous in $t$ and $x$
when restricted to the curve\cite{Li.thesis}. In this case  $T_xF_t(v)$ is the
 derivative  along a curve tangent to $v_0$ at $x_0$.

Let $f$ be a bounded function with bounded continuous first derivative,
 we can differentiate
 $P_tf$  to get
\begin{equation}\label{4}
d(P_tf)(v)=\E df(TF_t(v)), \hskip 18pt v\in T_xM
\end{equation}
given strong 1-completeness and assuming 
$\sup_{x\in K}\E|T_xF_t|^{1+\delta}$ finite for all compact sets $K$ and
some $\delta>0$. Note that given nonexplosion
 $\sup_{x\in K}\E|T_xF_t|^{1+\delta}$ is finite  if
$H_{1+\delta}$, as defined in section 2, is bounded above.  See \cite{flow}.

Alternatively   define a linear map  $\delta P_t$ on differential
 1-forms by
\begin{equation}\label{5}
\delta P_t(\phi)(x)(v)=\E \phi(TF_t(v))\chi_{t<\xi(x)}, \hskip 6pt  v\in T_xM,
\end{equation}
 which is formally a semigroup due to the  Markov property of the
solution flow. Furthermore if  (\ref{4})  holds, then
$$\left(\delta P_t\right)(df)(v)=\E df(TF_t(v))=d(P_tf)(v).$$
 This equality of $\delta P_t$ and $dP_t$  turns out to be a very interesting
property of the SDE and will be extensively used in   section 3 and 4.

\bigskip

 \noindent {\bf B.}
 For $p$  an integer and $K$ a subset of the manifold, define the 
$p^{th}$-moment exponents as follows:
\begin{equation}
\mu_K(p)=\limsup_{t\to \infty} \sup_{x\in K}
{1\over t}\log \E|T_xF_t|^p.
\end{equation}

The SDE is said to be {\it $p^{th}$-moment stable }if $\mu_{x}(p)<0$ for 
all $x$ in $M$, {\it strongly   p$^{th}$-moment stable} if $\mu_{K}(p)<0$
 for all  compact sets $K$. It is p$^{th}$-moment unstable if $\mu_{x}(p)>0$
for all $x$.

\bigskip

\noindent{\bf C. Main results:}
 We are mainly interested in non-degenerate stochastic differential equations.
 Recall that (1) is said to be a {\it Brownian system (with drift
 $Z$)} if it has generator $\half \Delta$ ( $+Z$ ) for  $\Delta$ the
 Laplacian.  It is called a {\it gradient Brownian system (with  drift $Z$)}
 if $X$ is given by an isometric immersion $j:  M\to  \R^m$, i.e.
for each $e\in \R^m$ and $x\in M$, $X(x)(e)$ is given by $\nabla <j(x), e>$.
The solution flow to the SDE is then a Brownian flow or gradient Brownian
 flow respectively. If $Z=\nabla h$ for some function $h$ on the manifold,
 then we  have {\it  h-Brownian systems}.

\bigskip

Let Ric$_x$ denote the Ricci curvature at $x$, and
\[  \begin{array}{ll}
h_p(x)=-p\, \inf_{ |v|\le 1}\left\{\right.
 &Ric_x(v,v)-2\left<\nabla Z(x)(v), v\right>
-\sum_1^m\left|\nabla X^i(v)\right|^2 \\
& -(p-2)\sum_1^m{\left<\nabla X^i(v), v\right>^2\over |v|^2}\left.\right\}
 \end{array} \]

We  say \cite{EL-LI-RO} that a function $f$ is {\it strongly stochastically
 positive} if for all  compact subsets $K$,
$$\limsup_{t\to \infty}{1\over t}\sup_{x\in K}\log \E
\left(\exp^{-\half \int_0^t f(F_s(x))ds}\right)<0.$$ 
It is {\it stochastically positive} if the above holds for all $K=\{x\}$ for
 all $x$ in $M$.

\bigskip

\noindent {\bf Theorem 2.2   }
 A complete gradient h-Brownian system is p$^{th}$-moment  stable if $-h_p$
 is  stochastically positive;  strongly p$^{th}$-moment  stable if
 $-h_p$ is strongly stochastically positive.

\bigskip

\noindent{\bf Theorem 4.1 }
The fundamental group $\pi_1(M)$ of a complete Riemannian manifold vanishes
if on it there is a strongly 1-complete strongly moment stable h-Brownian
 system with the property $d(P_tf)=\delta P_t(df)$ for all $f$ in $C_K^\infty$,
the space of smooth functions with compact support.

\bigskip

In particular,

\noindent{\bf Corollary 4.2 }
The first fundamental group of a complete Riemannian manifold  vanishes
if there is a h-Brownian motion on it such that $\sup_{x\in M}h_1(x)$
is negative.

\bigskip

In terms of the geometrical quantities of the manifold, we have: 
 let $r$  be the distance function between $x$ and a fixed point of the
 manifold:

\noindent{\bf Corollary 4.3 } Let $M$ be a closed submanifold of $\R^m$ with
 its second fundamental form $|\alpha_x|$ bounded by
 $c\left[1+\ln(1+r(x))\right]^\half$.
Let $h$ be a smooth function with  
$|{\partial h\over \partial r}|\le c[1+r(x)]$ and 
$Hess(h)\le  c\left[1+\ln(1+r(x))\right]$. 
Then $\pi_1(M)=\{0\}$ if $-h_1$ is strongly stochastically positive.

\bigskip

For analogous results when $M$ is compact, including the vanishing of
higher homotopy groups see \cite{EL-RO96}, where we also give an extension to tensor products
of the derivative flows.

 We also discuss the existence of  moment stable Brownian systems on $\R^n$:

\noindent{\bf Corollary \ref{R1}: } 
 There is  no moment stable Brownian system on $\R^1$.

\section{Conditions for strong moment stabilities}

\noindent {\bf A.} 
Let $\{e_1, e_2, \dots, e_m\}$ be an orthonormal basis for $\R^m$. Write
$X^i(x)=X(x)(e_i)$ and $B_t=(B_t^1, \dots, B_t^m)$ giving $m$ independent
 1-dimensional Brownian motions. The (1) can be written as:
$$dx_t =\sum_1^m X^i(x_t)\circ dB_t^i +A(x_t)dt$$
and (2) as
$$dv_t =\sum_1^m \nabla X^i(v_t)\circ dB_t^i +\nabla A(v_t)dt.$$

\noindent For   $v\in T_xM$,  define 

\begin{equation}\label{hp}\begin{array}{ll}
H_p(x)(v,v)&= 2\left<\nabla A(x)(v),v\right>+
\sum_1^m\left <\nabla^2 X^i(X^i,v),v\right>
+\sum_1^m \left<\nabla X^i\left(\nabla X^i(v)\right),v\right>\\
&+\sum_1^m \left|\nabla X^i(v)\right|^2
+(p-2)\sum_1^m{1\over|v|^2}\left<\nabla X^i(v), v\right>^2.
\end{array}\end{equation}

If the SDE is a Brownian system with drift $Z$, then  the first 3 terms of
$H_p$ becomes $-Ric_x(v,v)+2<\nabla Z(x)(v),v>$.

 By an It\^o's formula from \cite{ELflow},
\begin{eqnarray*}
|v_t|^p&=&|v_0|^p+p\int_0^t |v_s|^{p-2}\sum_1^m
\left <\nabla X^i(v_s), v_s\right> dB_s^i\\
&&+{p\over 2}\int_0^t |v_s|^{p-2}H_p(x_s)(v_s,v_s)ds.\\
\end{eqnarray*}

\noindent
Letting  
 $$M_t^p= \sum_1^m \int_0^t p  {<\nabla X^i(v_s), v_s>_{x_s} \over |v_s|^2} 
dB_s^i, $$ 

$$a_t^p={p\over 2} \int_0^t {H_p(x_s)(v_s,v_s) \over |v_s|^2} ds, $$
and ${\cal E} (M_t^p)=\exp^{M_t^p-{<M^p, M^p>_t \over 2}}$, as in \cite{flow}
 we solve the  equation for $|v_t|^p$  to get:

\begin{equation}\label{eq: vt as exponential}
|v_t|^p=|v_0|^p{\cal E}(M_t^p)\exp^{a_t^p}. \end{equation}
C.f. \cite{TANI89}.
\bigskip

Clearly if $H_p$ is negative, i.e. $H_p(x)(v,v)< -c^2|v|^2 $ for some
 $c\not =0$, then 
$$\limsup_{t \to \infty}{1\over t}\log \sup_{x\in M}\E |T_xF_t|^p<0$$
as known\cite{ELflow}. But for gradient Brownian systems, we can do  better.

\noindent {\bf B.}
Let $j: M\to \R^m$ be the isometric immersion giving rise to the gradient
 Brownian system and let $\nu_x$ be the space of normal vectors at $x$. 
 There  is the second fundamental form:
$\alpha_x: T_xM\times T_xM \to \nu_x$ and the shape operator:

$$A_x: T_xM\times \nu_x \to T_xM$$

\noindent
related by $\left<\alpha_x(v_1,v_2), w\right> =\left<A_x(v_1,w), v_2\right>$.
If $Y(x): R^m\to \nu_x$ is the orthogonal projection, then

$$\nabla X^i(v)=A_x \left(v,Y(x)e_i\right)$$

\noindent
as shown in \cite{ELbook} and \cite{ELflour}.
 Denote by $|\alpha_x(v,\cdot)|_{H,S}$ the corresponding Hilbert Schmidt norm,
 and $|\cdot|_{\nu_x}$ the norm  in $\nu_x$.  Then

$$\sum_1^m|\nabla X^i(v)|^2=|\alpha_x(v,-)|^2_{H,S},$$
and $$\sum_1^m<\nabla X^i(v),v>^2=|\alpha_x(v,v)|^2_{\nu_x},$$
for $v\in T_xM$. This gives\cite{ELflow}:

\begin{equation}\label{eq: gradient H}
\begin{array}{cl}
H_p(x)(v,v)=& -{\rm Ric}_x(v,v)+2<\nabla A(x)(v),v>+
|\alpha_x(v,\cdot)|_{H,S}^2\\
 &+ {(p-2) \over |v|^2}|\alpha_x(v,v)|^2_{\nu_x}.
\end{array}
\end{equation}

\bigskip

\noindent
Let $h_p(x)= \sup_{|v|\le 1}p\,H_p(x)(v,v)$, and 
$\underline h_p(x)=\inf_{|v|\le 1}p\,H_p(x)(v,v)$.

\begin{lemma} \label{le: moment stability}
 Consider a   gradient h-Brownian system. Assume nonexplosion. Then
$$\E\left(\exp^{\half \int_0^t \underline h_p(F_s(x))dx}\right)
\le \E|T_xF_t|^p\le n\E\left(\exp^{\half \int_0^t h_p(F_s(x))ds}\right).$$
\end{lemma}

\noindent {\bf Proof:}
Let $Y(x)^*: \nu_x \to \R^m$ be the adjoint of the orthogonal projection
 $Y(x)$, then $X(x)\left(Y(x)^*(\alpha _x(v,v)\right)\equiv 0$ and
\begin{eqnarray*}
&& <Y^*(\alpha_x(v,v)), dB_t>=<\alpha_x(v,v), Y(dB_t)>\\
=&&\sum_1^m <A_x(v, Y(e_i)), v>dB_t^i=\sum_1^m <\nabla X^i(v), v>dB_t^i
\end{eqnarray*}
from the  definitions. Let
$$\tilde {B_t} =B_t-\int_0^t p\, Y^*\left(
\alpha_x({v_s\over |v_s|},{v_s\over |v_s|})\right)ds.$$
 Consider SDE  (2) on TM with $\{B_t\}$ replaced by $\{\tilde B_t\}$:

$$d\tilde v_t=\delta X(\tilde v_t)\circ d\tilde B_t +\delta A(\tilde v_t)dt.$$

Then  $\tilde x_t$  is  the  projection of $\tilde v_t$ to $M$ and
  solves the stochastic differential equation:
\begin{equation}
d\tilde x_t=X(\tilde x_t)\circ d\tilde{B_t}+A(\tilde x_t)dt,
\label{eq:Girsanov}
\end{equation}
which is just   SDE (1) from 
 $X(x)\left(Y(x)^*(\alpha _x(v,v))\right)\equiv 0$. So  $\tilde x_t$ has
 the same  distribution as $x_t$.

 Then by (\ref{eq: vt as exponential}) and the Cameron-Martin-Girsanov formula,

\begin{eqnarray*}
\E|v_t|^p&=&|v_0|^p\E 
\exp^{M_t^p -\half<M^p, M^p>_t}
\exp^{ \int_0^t {p\,H_p(x_s)(v_s,v_s) \over 2|v_s|^2} ds}\\
&\le&  |v_0|^p\E  \exp^{M_t^p -\half<M^p, M^p>_t}
 \exp^{\int_0^t{h_p(x_s)\over 2} ds}\\
&=&|v_0|^p\E \exp^{\int_0^t {h_p(\tilde x_s)\over 2} ds}
= |v_0|^p\E \exp^{ \int_0^t {h_p(x_s)\over 2 } ds}.
\end{eqnarray*}

Thus
\begin{equation}\label{eq: Girsanov 25}
\E|T_xF_t|^p\le n\E \left(\exp^{\half \int_0^t h_p(F_s(x))ds}\right).
\end{equation}
The other half of the required inequality follows in the same way.
\hfill\rule{3mm}{3mm}

\begin{theorem}
 A gradient SDE  is $p^{th}$-moment stable if

\begin{equation}
\limsup_{t\to\infty}{1\over t}
\log \E\left(\exp^{\half \int_0^t h_p(F_s(x))ds}\right)<0
\end{equation}
for each $x$ in $M$, and it is $p^{th}$-moment unstable if
$$\liminf_{t\to \infty}{1\over t}\log
 \E\left(\exp^{\half \int_0^t\underline h_p(F_s(x))ds}\right)>0.$$
Similarly  it is strongly $p^{th}$-moment stable if
$$\limsup_{t\to\infty}{1\over t} \sup_{x\in K}\log
  \E\left(\exp^{\half \int_0^t h_p(F_s(x))ds}\right)<0$$
for each compact subset $K$.

\end{theorem}

\section{The Recurrency of h-Brownian Motions}

{\bf A.} Let $h: M\to \R$ be a smooth function on $M$ and
 $\Delta^h=\Delta+2L_{\nabla h}$ be the Bismut-Witten Laplacian for
 $L_{\nabla h}$ the Lie derivative in the direction of $\nabla h$. Then
it is known that $\Delta^h$ is an essentially self-adjoint operator.
Its closure shall also be denoted by $\Delta^h$. Moreover if $\heatsemi$
is the semigroup defined by the spectral theorem with $\heatsemif$
its restriction on differential 1-forms, then
 $d\left(\heatsemi f\right)=\heatsemif (df)$
for $f\in C_K^\infty$, the space of smooth functions with compact support.
A SDE with generator $\half \Delta^h$ is called an h-Brownian system, the
 solution  is  an h-Brownian motion.

As usual $P_t$ is the transition semigroup associated with our SDE, which
equals  the heat semigroup $\heatsemi$ on bounded $L^2$ functions
 \cite{El-RO88}.  We shall show that  a complete strongly moment stable
 h-Brownian motion  is recurrent  if
 $d(P_t f)=\left(\delta P_t \right) (df)$ holds on $C_K^\infty$.

\noindent {\bf B.}
 But the assumption on the   completeness of the h-Brownian   motion is not
 really an extra condition  since  $d(P_t  f)=\left(\delta P_t\right)(df)$
 for $f$ in $C_K^\infty$ implies nonexplosion if
 $\E|T_xF_t|\chi_{t<\xi(x)}<\infty$ on an open set $U$
and for $t<t_0$ for some constant $t_0$.
In fact \cite{Li.thesis} it is known that an h-Brownian motion on a complete
 Riemannian   manifold does not explode if  there is an open set $U$ and a
 number $t_0>0$   such that
 $$|\heatsemif df|_x\le c_{t}(x)|df|_\infty, \hskip 6pt x\in U, \hskip 3pt
 t\le t_0. $$
Here $c_t(x)$ is a constant depending possibly on $t$ and $x$ but not on the
 function $f$ in $C_K^\infty$.    But if $d(P_tf)=\delta P_t(df)$   we can 
obtain an upper  bound for the  $L^\infty$ norm of the heat  semigroup   on
 differential 1-forms in terms  of $|TF_t|$:
\begin{eqnarray*}
\left|\heatsemif (df)(v)\right|&=&\left|d(P_tf)(v)\right|=
\left|\E df(TF_t(v))\chi_{t<\xi}\right|\\
 &\le& |df|_\infty|v|\E\left|T_xF_t\right|\chi_{t<\xi} \end{eqnarray*}
and so take $c_t(x)=\E|T_xF_t|\chi_{t<\xi}$.

\bigskip

We should note here a more standard non-explosion criteria which can be
deduced by the same type of argument.
Using  another  expression of the heat semigroup  on 1-forms:
$\heatsemif(df)=\E df(W_t^h)\chi_{t<\xi}$ for $W_t^h$  the Hessian flow
 (see pp. 360-362 of  \cite{ELflour}), we conclude:
 an  h-Brownian motion does not explode if 
$$\E\sup_{t\le T}\exp^{\half\int_0^t\rho^h(x_s)ds}\chi_{t<\xi(x_0)}<\infty,
\hskip  10pt x\in U,$$
for  some open set $U$, constant $T>0$.  Here  
\begin{equation}
\rho^h(x)=
-\inf_{|v|= 1, v\in T_xM}  \left\{ Ric_x(v,v)-2Hess(h)(v,v)\right\}.
\end{equation}
This extends Bakry's result\cite{Bakry86}: an  h-Brownian motion does not
 explode if  $\rho^h$ is bounded  above.

\bigskip

\noindent{\bf C. }
 On the other hand given nonexplosion, there is an easy test for
 $d(P_tf)=\delta P_t(df)$  for Brownian systems with drift $Z$\cite{Li.thesis}.
 We only  need to check  $\E\sup_{s\le t}|T_xF_s|<\infty$.
This very condition on $ |T_xF_s|$ also gives strong 1-completeness and
can be realized by assuming $H_1(x)$ given by (\ref{hp}) is bounded above.
 Note that $H_p$ is  nondecreasing as $p$ increases, and for $p\ge 1$
$$H_p(x)(v,v)\ge -Ric_x(v,v)+2<\nabla Z(v)v,v>.$$
So for h-Brownian motions the condition that $H_p$ is bounded above for
 some $p\ge 1$ gives  nonexplosion,  and therefore strong 1-completeness and
 $d(P_tf)=\delta P_t(df)$ for all functions $f\in C_K^\infty$.

\bigskip

\noindent {\bf D.} Let $M$ be a complete Riemannian manifold and  $dx$ its
 Riemannian volume element. Then it is known that  $\exp^{2h} dx$ is   the
 invariant measure for the h-Brownian motion.  The manifold  is  said to
 have finite h-volume if   $\int_M \exp^{2h(x)}dx<\infty$.
 Arguing by contradiction we see that  the  h-Brownian motion  on $M$ is
 recurrent if  $M$ has finite  h-volume.

\begin{proposition}\label{finite volume}
 Suppose there is a   h-Brownian system
such that  $d( P_t f)=\delta P_t(df)$ on $C_K^\infty$ and 
$$\int_0^\infty \sup_{x\in K} \E|T_xF_t|\chi_{t<\xi(x)}dt<\infty.$$
 Then $M$ has finite  h-volume.
\end{proposition}

For the proof, we mimic  Bakry \cite{Bakry86}.
Let $\{h_n\}$ be a sequence of increasing functions bounded between 0 and 1
with limit 1 and $|\nabla h_n|\le {1\over n}$.  For $f\in C_K^\infty$, 
let $P_\infty f$
be the limit of $\heatsemi f$ as $t$ goes to infinity. Then $P_\infty f$  is
 harmonic  and thus a constant. Suppose the h-volume is not finite, then 
$P_\infty f$  has to be zero. We shall show that this is impossible.

Take $g\in C_K^\infty$, then
$$\int_M(P_\infty f-f)ge^{2h}\, dx
=\lim_{t\to \infty} \int_M(P_tf-f)ge^{2h}\, dx.$$

But 

\begin{eqnarray*}
&& \int_M (\heatsemi f-f)ge^{2h}dx 
=\half \int_0^t \int_M \Delta^{h}\left(\exp^{\half s\Delta^h}f\right)g
e^{2h}dx ds\\
  &&=\half \int_0^t\int_M< d(\exp^{\half s\Delta^h} f), dg> e^{2h}dxds
=\half \int_0^t\int_M < \exp^{\half  s\Delta^{h,1}}  (df),dg> e^{2h}dx ds\\
  &&= \half \int_0^t\int_M <df,\exp^{\half  s\Delta^{h,1}} (dg)>e^{2h} dx ds\\
&&=\half
 \int_0^t\int_{Supp(f)} <df, \exp^{\half  s\Delta^{h,1}}(dg)>e^{2h} dx ds,
\end{eqnarray*}

\noindent
since $\heatsemif$ is self-adjoint.  Here $Supp(f)$ denotes the support of $f$.
 Note that there is no explosion by the  previous argument.  
Replacing $\heatsemif(dg)$ by $\E dg(TF_t)$ in the above calculation, we get

\begin{eqnarray*}
&&|\int_M (P_tf-f)ge^{2h}dx|
\le \half |\nabla g|_\infty\,
\int_0^t  \sup_{x\in Supp(f)} \E\left(|T_xF_s|\right) ds 
 \int_M|\nabla f| e^{2h} dx\\
\end{eqnarray*}
Letting  $g=h_n$ and taking $t$ to infinity,

$$|\int_M -f \,e^{2h}dx|
\le\lim_{n\to \infty} {1\over 2n} |\nabla f|_{L^1}
\int_0^\infty \sup_{x\in Supp(f)} \E\left(|T_xF_s|\right) ds= 0.$$

\noindent
 Since we can  choose a function   $f\in C_K^\infty$ with
 $\int_M fe^{2h} dx\not = 0$, this gives a contradiction. \rule{3mm}{3mm}

\bigskip

\noindent{\bf E.}
  Note that  in Elworthy and Yor \cite{EL-YOR} it was showed that  for $M$
compact  moment stability is impossible for Brownian systems if $Ric_x\le 0$
 at all points. In fact,  by (\ref{eq: vt as exponential}),
 there is no moment stable SDE on a manifold if
$H_p(x)$ defined by (\ref{hp}) is non-negative and if 
$\E\exp^{\half \int_0^t|\nabla X(x_s)|^2ds}<\infty$ for all $t$.
However
\begin{eqnarray*}
H_p(x)(v,v)& \ge& 2<\nabla A(x)(v), v>\\
&+&\sum_1^m <\nabla^2 X^i(X^i, v), v> +
\sum_1^m <\nabla X^i(\nabla X^i(v), v>,
\end{eqnarray*}
   by (\ref{hp}) and the right hand side equals
 $2<\nabla Z(x)(v), v> -Ric_x(v,v)$ for
Brownian systems with drift $Z$.  For $R^n$ the condition 
$\E\exp^{\half \int_0^t|\nabla X(x_s)|^2ds}<\infty$ for all $t$
 is satisfied if   $X$ and $A$ have  linear growth and $|D X|$ and $|D A|$
 have sub-logarithmic growth. See  \cite{flow}. 

Furthermore a complete   gradient Brownian system (with drift) is not moment
 stable if  $H_p\ge 0$, by lemma 2.1. In particular,  c.f. \cite{ELflow}: 

\begin{corollary}\label{R1}
There is no moment stable Brownian system on $\R^1$.
\end{corollary}

When the Ricci curvature is not non-positive, we have an analogous result here:
 let $r(x)=d(x, x_0)$ the distance of $x$ from a fixed point $x_0$ of $M$,

 \begin{corollary}
Suppose $M$ is complete and non-compact  with
$${\rm Ric_x} > -{n \over n-1} {1\over r(x)^2}, \hskip 6pt r(x)>r_0,
 \hskip 3pt x\in M$$
for a constant $r_0>0$. Then a Brownian system on $M$ cannot  be strongly
 moment stable and have  $d(P_tf)=\delta P_t(df)$ for all $f\in C_K^\infty$.
\end{corollary}

 This is an application of the following  from  \cite{ch-gr-ta}:
 The volume of $M$ is infinite   for  noncompact manifolds with the above
 condition on the Ricci curvature.
However we do not have a good extension of the corollary to h-Brownian
 systems(unless in the trivial case when the function $h$ is bounded).

As a consequence there is no strongly moment stable Brownian system on
 $\R^n$ with
 $\sum_1^m |D X^i(x)(v)|^2-\sum_1^m {1\over |v|^2}<D X^i(x)(v), v>^2$
bounded,  for  this gives  $d(P_tf)=\delta P_t(df)$ for $f\in C_K^\infty$
 by part C of this section.

\section {Vanishing of $\pi_1(M)$}

As pointed out by  Elworthy \cite{EL-survey}, the stability of a stochastic 
flow  is directly  related to the topological properties of the underlying  
manifold. In particular there is the following theorem: A Brownian motion
 on a compact Riemannian manifold cannot be moment stable if the first 
fundamental group $\pi_1(M)$ of $M$ is not trivial. It is proved with the 
help of the following: On a compact manifold with non-vanishing $\pi_1(M)$
there is a loop of strictly positive  minimal length in its homotopy class. 
This is not true in general for noncompact manifold. However for $M$ not
 compact but with positive injective radius, a similar
argument shows \cite{Li.thesis} that $\pi_1(M)=\{0\}$  if  there is a
 strongly 1-complete  h-Brownian system such that for all compact subsets $K$, 
 $\lim_{t\to \infty}\sup_{x\in K} \E|T_xF_t|=0$. 
The key for a proof  for  general noncompact manifolds   is the recurrency of 
 h-Brownian motions.
Recall that an h-Brownian motion has a natural invariant measure
 $\exp^{2h}dx$. Furthermore if the invariant measure is finite and 
  $\mu$ is the normalized invariant measure,  we have the  recurrency:
$$\lim_{t\to \infty}P_t(\chi_K)(x)=\mu(K).$$

\begin{theorem}\label{th: homotopy}
Let $M$ be a complete Riemannian manifold. Suppose there is a  strongly 
1-complete  h-Brownian  system such that
	$$\int_0^\infty \sup_{x\in K} \E|T_xF_t| dt<\infty.$$
Then the first homotopy group $\pi_1(M)$\index{$\pi_1(M)$} vanishes  if
 $dP_tf=\delta P_t(df)$ for $f\in C_K^\infty$ (or more generally
 if $M$ has finite  h-volume).
\end{theorem}

\noindent{\bf Proof:}
Take $\sigma$ to be a $C^1$ loop parametrized by arc length. Then 
$F_t\circ\sigma$ is a $C^1$ loop homotopic to $\sigma$ by the strong 
1-completeness. Let $\ell(\sigma_t)$ denote the length of $F_t(\sigma)$
with $\ell_0=\ell(\sigma_0)$.

If we can  show $F_t\circ \sigma$ is contractible  to a point in $M$ with
 probability bigger than zero for some $t>0$, then the theorem is proved 
from the definition:  $\pi_1(M)=0$ if every continuous loop is contractible
 to one point. 

First we claim there is a sequence of numbers $\{t_j\}$ converging to 
infinity  such that:

\begin{equation}
 \E\ell(\sigma_{t_j}) \to 0. 
\label{eq:homo1}
\end{equation}
Since:
\begin{eqnarray*}
\int_0^\infty \E\ell(\sigma_t)dt 
 &\le&\int_0^\infty \E\left(\int_0^{\ell_0} |T_{\sigma(s)}F_t|ds\right)dt\\
&\le&\ell_0 \int_0^\infty \sup_s \E|T_{\sigma(s)}F_t|dt<\infty.  
\end{eqnarray*}

\noindent
So $\liminf_{t\to \infty} \E\ell(\sigma_{t}) =0$, giving $(\ref{eq:homo1})$.
Therefore $\ell(\sigma_{t_j})\to 0$ in probability.

Note that by proposition \ref{finite volume}
 we have finite  h-volume  if $dP_tf=(\delta P_t)(df)$ for
 $f\in C_K^\infty$.  
Let $\mu$ be the normalized invariant measure on $M$ for the process. Let
 $K$ be a compact set in $M$ containing the image set of the loop $\sigma$
 and  which has measure $\mu(K)>0$.  Let   $a>0$ be the infimum over 
$x\in K$ of the  injectivity radius at $x$.

 Then first by $(\ref{eq:homo1})$, there is a number $N$ such that for $j>N$, 

$$P\{\ell(\sigma_{t_j})\ge\half a\} < {\mu(K) \over 4}.$$

\noindent
And then by ergodicity, 

$$\lim_{t\to \infty} P\{F_t(x)\in K\}=\mu(K),$$
for   $x\in M$. Take a point $\tilde x$ in the image of the loop $\sigma$.
 There exists a number $N_1$ such that if $j>N_1$, then: 

$$P\{F_{t_j}(\tilde x)\in K\}  >{\mu(K) \over 2}.$$
Thus
\begin{eqnarray*}
&&P\{\ell(\sigma_{t_j})<\half a, F_{t_j}(\tilde x)\in K\}\\
&  & =P\{F_{t_j}(\tilde x) \in K\} 
-P\{F_{t_j}(\tilde x)\in K, \ell(\sigma_{t_j}) \ge\half a\}\\ 
&  &\ge P\{F_{t_j}(\tilde x) \in K\} -P\{\ell(\sigma_{t_j})\ge\half a\}\\
&  & > {\mu(K) \over 4}.
\end{eqnarray*}

But by the definition of injectivity radius, there is a coordinate chart
 containing a geodesic ball of radius ${a}$ around $F_{t_j}(\tilde x)$
 if  $F_{t_j}(\tilde x)$ belongs to $K$.  So the whole loop
 $F_{t_j}\circ\sigma$  is contained in the same chart and  thus contractible
 to one point  with probability  $>{\mu(k)\over 4}$.
\hfill \rule{3mm}{3mm}

It is clear from the proof that we do not need to assume either the SDS is a 
h-Brownian system or the equivalence of $dP_t$ and $\delta P_t$
  provided there is the recurrency.

\bigskip

 The following is known for arbitrary systems on compact manifolds
 \cite{ELflow}:

\begin{corollary}\label{th: group} 
 Let $M$ be a complete Riemannian manifold. Consider a 
h-Brownian system on $M$. Then  $\pi_1(M)=\{0\}$, if
$H_1(v,v)<-c^2|v|^2$ for some constant $c\not =0$.
\end{corollary}

\noindent {\bf Proof:} By part C. of section 3, the condition on
$H_1$ gives us the required strong 1-completeness and
 $dP_t=\delta P_t$ on $C_K^\infty$.

\begin{corollary} \label{co: group}
Let $M$ be a closed submanifold of $R^n$ with its second fundamental form 
$\alpha_x$ bounded by $c[1+\ln(1+r(x))]^\half$. Let $h$ be a smooth function
 on $M$ with $|{\partial h\over \partial r}|_x\le c|1+r(x)|$ and
$Hess(h)(x)\le c[1+\ln(1+r(x))]$.
 Then the first fundamental group vanishes if $-h_1$ is strongly stochastically
 positive.
\end{corollary}

\noindent{\bf Proof:}  Under these conditions, the gradient h-Brownian system
 on $M$ is  strongly complete and $d(P_tf)=\delta P_t(df)$ for all bounded 
functions  with bounded first derivative \cite{flow}. 

\bigskip

 When $M$ is compact, this gives
 $\pi_1(M)=\{0\}$ if  $\triangle^{h}+h_1<0$. See \cite{EL-RO94}.
Note in \cite{group}, it was shown that a compact manifold  $M$ has finite
 fundamental group under a weaker assumption:  $\triangle^h+\rho^h<0$.

\bigskip

\noindent {\bf Remark:}

 A SDE is strongly p-complete if its solution $\{F_t(x)\}$ is jointly
continuous in $t$ and $x$ a.s.  when restricted to any smooth singular
 simplices.  If on $M$ there is a strongly p-complete SDE which is also
 strongly p$^{th}$-moment stable, then all bounded closed p-forms are exact.
 Consequently the  natural map from the p$^{th}$ real cohomology with
 compact support to the p$^{th}$ real cohomology is trivial for such
 manifolds \cite{Li.thesis}.

\bigskip

\noindent{\bf Acknowledgment:} Some of the results here first appeared in my
 Ph.D. thesis. I would like to thank Professor D. Elworthy for valuable
 suggestions and helpful discussions.

References are updated in this version.
\end{document}